\documentclass[11pt]{article}
\usepackage{amstext,amssymb,amsmath,amsbsy}

\usepackage{hyperref}
\usepackage{amscd}
\usepackage{amsfonts}
\usepackage{indentfirst}
\usepackage{verbatim}
\usepackage{amsmath}
\usepackage{amsthm}
\usepackage{enumerate}
\usepackage{graphicx}
\usepackage[OT1]{fontenc}
\usepackage[latin1]{inputenc}
\usepackage[english]{babel}
\usepackage{amssymb}
\newtheorem{theorem}{Theorem}
\newtheorem{lemma}{Lemma}

\newtheorem{proposition}{Proposition}
\newtheorem{corollary}{Corollary}
\newtheorem{remark}{Remark}

\numberwithin{equation}{section}

\newcommand{\proofend}{\hfill $\Box$ }

\newcommand{\lr}{\rightarrow}
\newcommand{\dsp}{\displaystyle}

\newcommand{\supp}{\operatorname{supp}}
\newcommand{\dist}{\operatorname{dist}}

\newcommand{\dive}{\operatorname{div}}

\newcommand{\eps}{\varepsilon}

\newcommand{\loc}{_{loc}}

\newcommand{\mR}{\mathbb{R}}

\newcommand{\mc}{\mathrm{c}}

\numberwithin{equation}{section}

\title{On a regularized scheme for approximate acoustic cloaking using transformation optics} 

\author{Hoai-Minh Nguyen \footnote{School of Mathematics, University of Minnesota, MN, 55455,
hmnguyen@math.umn.edu} \footnote{The research is supported by NSF grant DMS-1201370 and by the Alfred P. Sloan Foundation.}}

\begin{document}

\maketitle

\begin{abstract} 

In this paper, we study approximate cloaking for the Helmholtz equation in $3d$  where the cloaking device is based on  transformations which blow up a cylinder of fixed height and small cross section of radius $\eps$ into the cloaked region. Assuming the zero Dirichlet boundary condition is imposed on the boundary of the cloaked region, we show that the degree of visibility is of order $\eps$ as $\eps$ goes to 0. This fact is quite surprising since 
it is known that  the degree of visibility, for  the scheme using transformations which blow up a small region of diameter $\eps$ into the cloaked region,  is of  order $\eps$ in $3d$ and $1/ |\ln \eps|$ in $2d$.  To understand the relation between these contexts, we as well revisit the known estimates and show that 
 the degree of visibility is of  order $\eps^{d-1}$ ($d=2,3$) for the scheme using  transformations which blow up a small {\it ball} of diameter $\eps$ into the cloaked region as long as the zero Dirichlet boundary condition is imposed on the boundary of the cloaked region. The symmetry of the cross section and of the ball are crucial so that these results hold. 

\end{abstract}  

\section{Introduction}

Cloaking via change of variables was introduced by Pendry, Schurig, and Smith \cite{PendrySchurigSmith} for the Maxwell system, and Leonhardt \cite{Leonhardt} in the geometric optics setting. They used a singular change of variables which blows up a point into a cloaked region.
The same transformation had been used before by Greenleaf, Lassas, Uhlmann \cite{GreenleafLassasUhlmann} to establish the non-uniqueness of Calderon's problem.  This singular structure implies not only difficulties in practice, but also in analysis. To avoid using the singular structure, regularized schemes have been proposed in \cite{Schurig06,Kildishev07,YanRuanQiu07,RuanYanNeffQiu,KohnShenVogeliusWeinstein,GreenleafKurylevLassasUhlmann07SHS}.

\medskip
Let us briefly review some known facts about approximate cloaking for the Helmholtz equation for a finite range, and for the full range of frequencies based on transformations which blow up a small ball of radius $\eps$ into the cloaked region.  This scheme was suggested by Kohn, Shen, Vogelius, and Weinstein in  \cite{KohnShenVogeliusWeinstein}. 
In this paragraph, whenever approximate cloaking is achieved, the degree of visibility, established by the authors,  is of  order $\eps$ in $3d$ and $1/ |\ln \eps|$ in $2d$. 
In order to successfully achieve approximate cloaking it is often advantageous to introduce a lossy layer in addition to the standard (mapped) cloak.
Using an appropriate lossy layer, it is proved in   \cite{KohnOnofreiVogeliusWeinstein} that approximate cloaking works well on a bounded domain
regardless of the contents of the cloaked region. In \cite{Liu}, the author established approximate cloaking for the exterior problems imposing the zero Dirichlet boundary condition. In \cite{NguyenHelmholtz}, the author showed  that approximate cloaking works well in the whole space regardless of the contents of the cloaked region using a fixed lossy-layer (the result in \cite{Liu} can be derived from \cite[Lemma 2.2]{NguyenHelmholtz}). In both \cite{KohnOnofreiVogeliusWeinstein} and \cite{NguyenHelmholtz} the authors demonstrated the necessity of the lossy layer, in order to obtain
a degree of approximate cloaking (near-invisibility) that  is independent of the contents of the cloaked region.
The paper \cite{NguyenVogelius2} establishes precise estimates for the degree of near-invisibility at all frequencies, where the dependence on frequency is explicit.
These estimates are sharp and independent of the contents of the cloaked region. To be a little more specific: in the high frequency case, the lossy scheme works as well as in the finite frequency case. However, the estimates degenerate as frequency tends to 0. This follows from (or can be explained by) the fact
that the lossy effect becomes weaker and weaker, as frequency tends to 0. 
These results play an important role in obtaining the degree of visibility for approximate cloaking for the  wave equations in \cite{NguyenVogelius3}. 
Without a lossy layer the situation becomes quite complicated, as explored in \cite{NguyenPerfectCloaking}. For example, in the $3d$ non-resonant case, i.e., when $k^2$ is not an eigenvalue of the Neumann problem inside the cloaked region
(here $k$ denotes the wave number), the approximate scheme works well: cloaking is achieved (as the parameter of regularization $\eps$ goes to zero) and the limiting field inside the
cloaked region is the corresponding  solution to the homogeneous Neumann problem. In the $3d$ resonant case, the situation changes completely. Sometimes cloaking is achieved (for example when there is no source inside the cloaked region); nevertheless,
the limiting field inside the  cloaked region depends on the solution in the free space (cloaking but not shielding; see also \cite{GreenleafKurylevLassasUhlmann-PE}  for a similar situation). Sometimes cloaking is not achieved, and the energy inside the cloaked region tends to infinity as the parameter of regularization tends to 0. In the $2d$ non-resonant case, the limiting field inside the cloaked region inherits a non-local structure (see also \cite{LassasZhou}). In the $2d$ resonant case, cloaking sometimes is not achieved, and the energy inside the cloaked region can go to infinity.  These ``lossless'' facts are somewhat different from what is frequently asserted in the literature, namely that $(a)$ in $3d$, cloaking is always achieved,  the limiting field inside and outside the cloaked region are independent, and the energy of the field inside the cloaked region remains bounded, and $(b)$ in $2d$, the limiting field inside the cloaked region satisfies the corresponding Neumann problem. 
In \cite{GreenleafKurylevLassasUhlmann07,WederRigorous,WederRigousTimeDomain,LiuZhou}, the singular cloaking scheme was studied using the notion of {\it finite energy} solution. The improvement of degree of invisibility have been recently investigated in \cite{AmmariKHL-Conductivity,AmmariKHL-Helmholtz,LiuSun}. 

\medskip
Recently, quantum cloaking has been investigated quite extensively; see  \cite{GreenleafKurylevLassasUhlmann-QP, GreenleafKurylevLassasUhlmann-Q, GreenleafKurylevLassasLeonhardtUhlmann} and references therein.  We note here that these results for approximating cloaking for the Helmholtz equations mentioned above can be used to obtain similar ones for approximate quantum cloaking; thanks to the gauge transformation (see e.g. \cite[Section 5]{GreenleafKurylevLassasUhlmann-Q}). The trapped states introduced by Greenleaf et. al. in \cite{GreenleafKurylevLassasUhlmann-Q}  in quantum cloaking corresponds to the resonant case for the Helmholtz equation in 3d.  We emphasize that one can do quantum cloaking for scattered waves even in a trapped state for which  there exist almost trapped states \cite{GreenleafKurylevLassasUhlmann-QP}. It is a consequence of the fact that one can do cloaking for the Helmholtz equation for the 3d resonant cases when there is no source inside the cloaked region and cloaking device \cite[Theorem 1.4]{NguyenPerfectCloaking}. These do not contradict to each other. Indeed, there is the instability of approximate cloaking. Nevertheless, the instability of cloaking takes place only if the potential varies as a function of parameter of regularization (even for a very small amount;  see  \cite{GreenleafKurylevLassasUhlmann-QP}).  The same phonemenon holds for the acoustic setting (see \cite[Proposition 1.13 and Remark 1.14]{NguyenPerfectCloaking}; in \cite{KohnOnofreiVogeliusWeinstein}, Kohn et. al.  previously observed the similar facts for bounded domain).  


\medskip
In this paper, we investigate approximate cloaking for the Helmholtz equation based on transformations which blow up a cylinder of fixed height and  small cross section into the cloaked region in $3d$. This scheme is a natural regularization of the one blowing up an interval into the cloaked region discussed in \cite{LeonhardtTyc}. This regularization is in the spirit of the one  introduced in \cite{KohnShenVogeliusWeinstein}.  The advantage of this scheme in comparison with the one blowing up a small ball is that the (material) transformation used for the cloaking device is less singular, hence facilitates the construction of such devices. 
This is the main motivation of our work  and to our knowledge, it is the first time that such a regularized scheme is considered and analyzed. In this paper, we impose the zero Dirichlet boundary condition on the boundary of the cloaked region. This   boundary condition is a good approximation if a highly conducting layer is used between the cloaked region and the layer from the blow-up scheme; see e.g. \cite{NguyenVogelius2, NguyenPerfectCloaking} for a discussion. The analysis on different boundary conditions (including the Dirichlet boundary condition) for a highly conducting obstacle can be found in, e.g., \cite{HaddarJolyNguyen1}.  Let $\eps$ be the parameter of regularization i.e., $\eps$ is radius of the small cross section. We show that the degree of visibility is of order of $\eps$.  We recall here that the degree of visibility in the 3d case is also of order $\eps$ if one uses the scheme which blows up a small region  of diameter $\eps$ into the cloaked region as mentioned above;  Clearly, this scheme is more singular than the one which blows up small cylinders. An important fact which makes our result hold is the symmetry of the cross section of the cylinder.  To understand the relation between these contexts, we also revisit the known estimates and show that  the degree of visibility is of order $\eps^{d-1}$ ($d=2,3$) when one uses transformations which blow up a small {\it ball} of diameter $\eps$ into the cloaked region as long as the zero Dirichlet boundary condition is imposed on the boundary of the cloaked region. Once again, the symmetry of the ball is crucial. 

\medskip
One might wonder what happens if one uses the ``regularized'' schemes which blow up a thin box of fixed length and width into the cloaked region in $3d$ or which blow up a small rectangle of fixed length into the cloaked region in $2d$. These are natural regularizations for the singular schemes which blow up a rectangle into the cloaked region in $3d$ and an interval into the cloaked region in $2d$. As long as,  the zero Dirichlet boundary is imposed on the boundary of the cloaked region, one can {\it not} use these schemes to do cloaking  since the trace of the solutions in the free space on these sets are not 0 in general and these sets are 1-codimentional -
not negligible in the trace sense.


\medskip
Let us describe briefly the idea of the proofs. Using transformation optics rule (or the change of variables formula - see e.g. Proposition~\ref{pro-TO}), it suffices to show that: 
\begin{itemize}
\item[i)] the effect of a cylinder of fixed height and small cross section of radius $\eps$ is of order $\eps$ in $3d$; 

\item[ii)] the effect of a small {\it ball} of radius $\eps$ is  of  order $\eps^{2}$ in $3d$, and $\eps$ in $2d$.  
\end{itemize}

The effect of small inclusions has been investigated intensively in the literature see, e.g.,  \cite{FriedmanVogelius,Cedio-FengyaMoskowVogelius,AmmariKang-Po,CapdeboscqVogelius,NguyenVogelius,Capdeboscq}. However, the way to obtain i) is new; the usual blow-up technique  does not work. The proof is processed as follows. 
 Considering the difference of the solutions in the free space and in the exterior of the small cylinder, we note that it is a solution to the homogeneous Helmholtz equations outside the cylinder and its value on the boundary of the cylinder is well-controlled. Hence it is natural to search bounds for its normal derivative so that the representation formula for the Helmholtz equation can be used. For this end, we use the Morawetz multiplier (Lemma~\ref{lem-Morawetz} in Section~\ref{sect-Pro1}; noting that the cylinder is star-shaped) and derive 
that the normal derivative behaves like, roughly speaking,  $1/ \eps$ on the boundary (see Lemma~\ref{lem-est-L2} and \eqref{est1-k} in Section~\ref{sect-Pro1}). Such an estimate is sufficient to obtain $\eps$ - degree of visibility from the contribution of the top and the bottom of the cylinder since their area is of order $\eps^{2}$. However, it is not good enough to reach $\eps$ - degree of visibility from the contribution of the rest of the boundary since its area is of order $\eps$ (the representation formula only gives an estimate of  order 1 for the difference). Here the symmetry of the cylinder plays a role. We observe that  the difference of the solutions is almost constant on each cross section since the cross section is small. By symmetry, it follows that the average of its normal derivative on each cross section is almost 0. 
This can be obtained from~Lemmas~\ref{lem-sym} (on the symmetry) and \ref{lem-stability-k} (on the stability) in Section~\ref{sect-Pro1}; Lemma~\ref{lem-stability-k} can be seen as a consequence of the fact that the singularity of solutions to the (Laplace) Helmholz equation on an interval is removable in $3d$. Taking into account this cancelation, by a standard averaging argument, statement i) follows.  Statement ii) can be obtained (similarly) using the representation formula, the standard estimates for small inclusions, and the fact that the average of the normal derivative on the boundary of the ball approximately equals $0$; this fact is a consequence of the symmetry of the ball as in the cylinder case mentioned above. 

\medskip 
We note that there are also other methods to do cloaking, e.g.,  using negative index materials,  plasmonic shells, 
or active cloaking devices. Concerning negative index materials,  one approach is based on the anomalous localized resonance introduced by Milton and Nicorovice in \cite{MiltonNicorovici} (see \cite{ MiltonNicoroviciMcPhedranPodolskiy,Milton-folded,BouchitteSchweizer10, AmmariCiraoloKangLeeMilton} references therein for further studies);  one approach is based on the concept of complementary media introduced by Lai et. al. in \cite{LaiChenZhangChanComplementary} (see \cite{NguyenNegative} for a more general setting and the  analysis). Cloaking using plasmonic shells is introduced by Alu and Engheta in \cite{AluEngheta} (see \cite{AluEngheta2,Tricarico-Alu}  and references therein for recent development). Cloaking using active devices is introduced by Miller in \cite{Miller} and  Guevara Vasquez et. al. in \cite{VasquezMiltonOnofrei} (see also \cite{VasquezMiltonOnofrei2}). 
 
\medskip
We now state our results precisely. For $\eps > 0$, let $C_{\eps}$ denote the set
\begin{equation*}
C_{\eps} : = \{(x', z) \in \mR^2 \times \mR;\; |x'| < \eps \mbox{ and }     z \le 1/2\},
\end{equation*}
Let $D \subset \subset B_{2}$ be a open smooth  subset of $\mR^{3}$ containing the origin such that $\mR^{3} \setminus D$ is connected. 
Here and in what follows in $\mR^{d}$ ($d=2,\, 3$), $B_{r}$ ($r >0$) denotes the ball centered at the origin of radius $r$.
For $\eps > 0$, let $F_{\eps}$ be a bi-Lipschitz transformation \footnote{This means $F$ is bijective, and $F$ and $F^{-1}$ are Lipschitz.} from $\mR^{3}$ onto $\mR^{3}$ such that 
\begin{equation*}
F_{\eps} (C_{\eps}) = D, \quad F_{\eps}(\partial C_{\eps}) = \partial D,  \quad \mbox{ and } \quad F_{\eps}(x) = x \mbox{ if } |x| \ge 2. 
\end{equation*} 

Here is an example. Assume $D= D_{1} = \{(x', z) \in \mR^{2} \times \mR; \; |x'| \le 1/2; |z| \le 3/4\}$ and let $D_{2} = 2 D_{1} = \{(x', z) \in \mR^{2} \times \mR; \; |x'| \le 1; |z| \le 3/2 \}$.  
Then one can choose $F_{\eps}$ as follows: 
\begin{equation}\label{defFeps}
F_\eps(x', z) = \left\{\begin{array}{cl} (x',z) & \mbox{ if } (x',z) \in \mR^3 \setminus D_{2}~, \\[12pt]
\dsp \Big(\Big[\frac{1 - 2\eps}{2(1 - \eps)} + \frac{|x'|}{2(1 - \eps)} \Big] \frac{x'}{|x'|},  \frac{3}{4}z + \frac{3}{8} \Big) & \mbox{ if } (x', z) \in D_2 \setminus C_{\eps}, \\[12pt]
\dsp \Big(\frac{x'}{2\eps}, \frac{3}{2}z \Big)  & \mbox{ if } (x', z) \in C_\eps.
\end{array}\right.
\end{equation}
 
Define 
\begin{equation*}
A_{c}, \Sigma_{c} = \left\{ \begin{array}{cl} 
I, 1 & \mbox{ for } |x| \ge 2, \\[6pt]
{F_{\eps}}_{*} I, {F_{\eps}}_{*} 1 & \mbox{ for } x \in B_{2} \setminus D.
\end{array} \right.
\end{equation*}
Here and in what follows we  use the standard notation
\begin{equation}\label{defF*}
F_*A(y) = \frac{D F (x) A(x) DF^T(x)}{| \det DF(x) |} \quad \mbox{and} \quad  F_*\Sigma(y) = \frac{\Sigma(x)}{ |\det DF(x)|}, \quad \mbox{ with } x = F^{-1}(y)~,
\end{equation}
for any symmetric matrix-valued function $A$, and any complex function $\Sigma$.

\medskip
Let $u \in H^{1}_{\loc}(\mR^{3})$ and $u_c \in H^{1}_{\loc}(\mR^{3} \setminus D)$ be the unique outgoing solutions to 
\begin{equation}\label{eq-u}
\Delta u + k^2 u = f \mbox{ in } \mR^3,
\end{equation}
and
\begin{equation}\label{eq-uc}\left\{
\begin{array}{cl}
\dive(A_c(x) \nabla u_c) + k^2 \Sigma_c (x) u_c = f & \mbox{ in } \mR^3 \setminus D,\\[6pt]
u_c = 0 & \mbox{ on } \partial D.
\end{array}\right.
\end{equation}
We recall that a function $v$ satisfies the well-known outgoing condition in $\mR^{d}$ $(d=2, \, 3)$ with respect to frequency $k$, means 
\begin{equation*}
\dsp \frac{\partial v}{\partial r} - i k v= o(r^{-(d-1)/2}) \quad \mbox{ as } r = |x| \lr \infty.
\end{equation*}

\medskip
Here is the main result of this paper.

\begin{theorem}\label{thm1} Let $k > 0$, $\eps > 0$, and $f \in L^2(\mR^3)$ be such that $\supp f \subset B_3 \setminus B_2$. Assume that   $u \in H^{1}_{\loc}(\mR^{3})$ and  $u_{c} \in H^{1}_{\loc}(\mR^{3} \setminus D)$ are the unique outgoing solutions to \eqref{eq-u} and \eqref{eq-uc}. We have
\begin{equation*}
\| u_c  - u \|_{H^1(B_5 \setminus B_2)} \le C \eps \|f \|_{L^2},
\end{equation*}
for some positive constant $C$ depending only on $k$. 
\end{theorem}

From the property of transformation optics (the change of variables formula, Proposition~\ref{pro-TO}), it follows that $u_{c} = u_{\eps}$ for $|x| > 2$, where $u_{\eps}$ is the solution to \eqref{eq-cl1} below.  Hence Theorem~\ref{thm1} is a consequence of  the following proposition which gives an estimate on the effect of the small inclusion $C_{\eps}$.

\begin{proposition}\label{pro1}
Let $k>0$, $0 < \eps < 1$,  $f \in L^{2}(\mR^{3})$ with $\supp f \subset B_{3} \setminus B_{2}$. Assume that   $u_{\eps} \in H^{1}_{\loc}(\mR^{3} \setminus C_{\eps})$ is the outgoing solution to the equation 
\begin{equation}\label{eq-cl1}
\Delta u_{\eps}  + k^{2} u_{\eps} =f \mbox{ in } \mR^{3} \setminus C_{\eps} \quad \mbox{ with } \quad u_{\eps} = 0 \mbox{ on } \partial C_{\eps}. 
\end{equation}
We have 
\begin{equation*}
\| u_{\eps} - u\|_{H^{1}(B_{5} \setminus B_{2})} \le C \eps \| f\|_{L^{2}},
\end{equation*}
for some positive constant $C$ depending only on $k$. Here $u \in H^{1}_{\loc}(\mR^{3})$ is the unique outgoing solution to \eqref{eq-u}. 

\end{proposition}

As in Theorem~\ref{thm1}, we can obtain the degree of visibility, which is of order $\eps^{d-1}$ ($d=2,3$), for the scheme blowing up a small ball of radius $\eps$ into the cloaked region. This fact is   a consequence of the following proposition:

\begin{proposition}\label{pro2}
Let $k>0$, $0 < \eps < 1$,  $f \in L^{2}(\mR^{d})$ ($d=2, \, 3$) with $\supp f \subset B_{3} \setminus B_{2}$. Assume that  $u_{\eps} \in H^{1}_{\loc}(\mR^{d} \setminus B_{\eps})$ is the unique outgoing solution to the equation 
\begin{equation}\label{eq-cl2}
\Delta u_{\eps}  + k^{2} u_{\eps} =f \mbox{ in } \mR^{d} \setminus B_{\eps} \quad \mbox{ with } \quad u_{\eps} = 0 \mbox{ on } \partial B_{\eps}. 
\end{equation}
We have 
\begin{equation*}
\| u_{\eps} - u\|_{H^{1}(B_{5} \setminus B_{2})} \le C \eps^{d-1} \| f\|_{L^{2}},
\end{equation*}
for some positive constant $C$ independent of $\eps$ and $f$. 
\end{proposition}

\medskip
To obtain Theorem~\ref{thm1} from Proposition~\ref{pro1}, we use the following standard proposition whose 
the proof is based on the change of variables formula and can be found, e.g.,  \cite[Section 2.2]{KohnOnofreiVogeliusWeinstein} or \cite[Proof of Proposition 1.1]{NguyenHelmholtz}.

\begin{proposition}\label{pro-TO} Let $k \ge 0$, $d=2, \, 3$, $D_1$ and $D_2$ be two open subsets of $\mR^d$, $F$ be a bi-Lipschitz from $D_1$ onto $D_2$, $a \in [L^\infty(D_1)]^{d \times d}$ and $\sigma \in L^\infty(D_1)$. Fix $u \in H^1(D_1)$ and set $v = u \circ T^{-1}$. We have
\begin{equation*}
\dive (a \nabla u) + k^2 \sigma u = f \mbox{ in } D_1,
\end{equation*}
if and only if
\begin{equation*}
\dive (F_*a \nabla v) + k^2 F_*\sigma v = F_* f \mbox{ in } D_2.
\end{equation*}
\end{proposition}

The paper is organized as follows. Sections~\ref{sect-Pro1} and \ref{sect-Pro2} are devoted to the proofs of Proposition~\ref{pro1} and Proposition~\ref{pro2} respectively. 


\section{Proof of Proposition~\ref{pro1}}\label{sect-Pro1}

\subsection{Preliminaries}

In this section, we will establish several lemmas  on the effect of the small inclusion $C_{\eps}$, which will be used in the proof of Proposition~\ref{pro1}. The first one concerns the average of the normal derivative of solutions on each cross section of $C_{\eps}$.

\begin{lemma}\label{lem-sym} Let $k, \eps > 0$, $g \in H^{1/2}(\partial C_{\eps})$,  and  let $v_\eps \in H^{1}_{\loc}(\mR^3 \setminus C_\eps)$ be the unique  outgoing solution to the equation
\begin{equation*}
\Delta v_\eps + k^{2} v_{\eps} = 0 \mbox{ in } \mR^3 \setminus \bar C_\eps,
\end{equation*}
with
\begin{equation*}
v_\eps = g \mbox{ on } \partial  C_\eps.
\end{equation*}
Assume that 
\begin{equation*}
g(x', z) = g(y', z) \quad \mbox{ for all } (x', z), (y', z) \in  (\mR^{2} \times \mR) \cap \partial C_{\eps}. 
\end{equation*}
We have
\begin{equation*}
\int_{\partial C_\eps ; z = a } \frac{\partial v_\eps}{\partial \eta} = 0 \mbox{ for } a.e. \, a \in (-1/2, 1/2).
\end{equation*}
\end{lemma}

\noindent{\bf Proof.} The conclusion is a consequence of the fact that
\begin{equation*}
v_\eps(x', z) = v_\eps(-x', z). 
\end{equation*}
This follows from the uniqueness of the outgoing solutions imposed the Dirichlet boundary condition.  
\proofend

\begin{remark}
In Lemma~\ref{lem-sym}, we only use the fact that the cross section of $C_{\eps}$ is radially symmetric. The constancy of the cross section is not necessary.  
\end{remark}

\medskip

The second lemma, whose the proof makes use of the Morawetz multiplier \cite{MorawetzLudwig} (see also \cite{Rellich} and \cite{PayneWeinberge}), will be applied to obtain estimates for the normal derivative of solutions on $\partial C_{\eps}$ in Lemma~\ref{lem-est-L2}. 

\begin{lemma} \label{lem-Morawetz} Let $k \ge 0$, $\Omega$ be a Lipschitz bounded open  subset of $\mR^{3}$ and 
$v \in C^{\infty}(\bar \Omega)$ be a solution to the equation $\Delta v + k^{2} v = 0$. 
We have 
\begin{equation*}
\frac{1}{2} \int_{\Omega} \big( |\nabla v|^{2} + k^{2} |v|^{2} \big) -   \int_{\partial \Omega} \Re \Big(  \langle \nabla v , \eta \rangle \langle \nabla \bar v , x \rangle \Big) =   \int_{\partial \Omega} \Re \Big( \frac{\partial v}{\partial \eta} \bar v + \frac{k^{2}}{2} \langle x, \eta \rangle |v|^{2} -  \frac{1}{2} \langle x, \eta \rangle |\nabla v|^{2} \Big). 
\end{equation*}
\end{lemma}

Here and in what follows $\langle \cdot, \cdot \rangle$ denotes the scalar product in $\mR^{3}$ and $\Re (a)$ denotes the real part of a complex number $a$.  In this lemma and its proof,  $\eta$ denotes the exterior unit normal vector on $\partial \Omega$. 

\medskip

\noindent
{\bf Proof.} In what follows, we only consider the real part of expressions. 
Multiplying the equation $\Delta v + k^{2} v = 0$ by $\bar v + x \cdot \nabla \bar v$ and integrating in $\Omega$, we obtain
\begin{equation}\label{t1}
0= -\int_{\Omega} |\nabla v|^{2} + \int_{\partial \Omega} \frac{\partial v}{\partial \eta} \bar v - \int_{\Omega} \nabla v \nabla \big[ x \cdot \nabla \bar v\big] 
+ \int_{\partial \Omega} \langle \nabla v, \eta \rangle \langle \nabla \bar v, x \rangle  + k^{2} \int_{\Omega} |v|^{2} 
+ k^{2} \int_{\Omega} v \big[ x \cdot \nabla \bar v\big] . 
\end{equation}
We have
\begin{equation}\label{t2}
k^{2} \int_{\Omega}  v \big[ x \cdot \nabla \bar v\big]  = \frac{k^{2}}{2} \int_{\Omega} x \cdot \nabla |v|^{2}  = -\frac{3k^{2}}{2} \int_{\Omega} |v|^{2} + \frac{k^{2}}{2} \int_{\partial \Omega} \langle x , \eta \rangle |v|^{2}. 
\end{equation}
Since $\nabla \big[ x \cdot \nabla \bar v\big] = \nabla \bar v + (x \cdot \nabla ) \nabla \bar v$, it follows that 
\begin{equation}\label{t3}
- \int_{\Omega} \nabla v \nabla \big[ x \cdot \nabla \bar v \big] =   - \int_{\Omega} |\nabla v|^{2} -\frac{1}{2} x \cdot \nabla |\nabla v|^{2} 
=  \frac{1}{2} \int_{\Omega} |\nabla v|^{2} - \frac{1}{2} \int_{\partial \Omega} \langle x, \eta \rangle |\nabla v|^{2}. 
\end{equation} 
A combination of \eqref{t1}, \eqref{t2}, and \eqref{t3} yields
\begin{multline*}
- \frac{1}{2} \int_{\Omega} \Big( |\nabla v|^{2} + k^{2} |v|^{2} \Big)  - \frac{1}{2} \int_{\partial \Omega} \langle x, \eta \rangle |\nabla v|^{2} \\[6pt]
+ \int_{\partial \Omega} \frac{\partial v}{\partial \eta} \bar v +  \frac{k^{2}}{2} \int_{\partial \Omega} \langle x , \eta \rangle |v|^{2} + \int_{\partial \Omega} \langle \nabla v, \eta \rangle \langle \nabla \bar v, x \rangle = 0 . 
\end{multline*}
Therefore, the conclusion follows. \proofend

\medskip

\medskip
The following lemma plays an important role in our analysis.

\begin{lemma}\label{lem-est-L2} Let $0 < \eps < 1$,  $g \in H^1(\partial C_\eps)$, and $v \in H^1_{\loc}(\mR^{3} \setminus \bar C_\eps)$ be the unique outgoing solution to the equation
\begin{equation*}
\Delta v + k^{2} v= 0 \mbox{ in } \mR^{3} \setminus \bar C_\eps \quad \mbox{ and } \quad v = g \mbox{ on } \partial C_\eps .
\end{equation*}
Then $\partial_\eta v \in L^2(\partial C_\eps)$. Moreover,
\begin{multline}\label{fund-est-L2}
C \left\{\int_{B_{4} \setminus C_{\eps}} \big( |\nabla v|^{2} + k^{2} |v|^{2} \big) +  \int_{\partial C_{\eps}}  \langle x, \eta \rangle |\partial_{\eta}  v|^{2} \right\}\\[6pt]
\le    \int_{\partial C_{\eps}} \frac{|x|^{2}}{\langle x, \eta \rangle} |\nabla_{T} g|^{2} + \int_{\partial C_{\eps}}  \frac{1}{\langle x, \eta \rangle} |g|^{2} + 
\| \nabla v\|_{L^{2}(\partial B_{4})}^{2} +  \| v\|_{L^{2}(\partial B_{4})}^{2}, 
\end{multline}
for some positive constant $C$ independent of $g$ and $\eps$. 
\end{lemma}

Here and in what follows, $\partial u/ \partial \eta$ and $\nabla_{T} u$ denotes the normal derivative and the tangential derivative of a function $u$ on $\partial D$, the boundary of a Lipschitz bounded  open subset of $\mR^{3}$,  whenever they are well-defined, and $\eta$ denotes the normal unit vector on $\partial D$ directed to the exterior of $D$ \footnote{We change the convention on $\eta$ from Lemma~\ref{lem-Morawetz}.} . 

\medskip

\noindent {\bf Proof.} The proof is given in the smooth case:  $v \in C^{\infty}(\overline{B_{4} \setminus C_{\eps}})$. The general case follows by approximation. In fact, by the same proof and the regularity theory of elliptic equations, one can show that the same result holds for $C_{\eps, \delta}$ and $g \in C^{\infty}(\partial  C_{\eps, \delta})$ where $C_{\eps, \delta}$ is a smooth star-shaped domain.  By choosing $C_{\eps, \delta}$ such that $C_{\eps, \delta}$ appropriately converges to $C_{\eps}$ as $\delta \to 0$, one obtains the conclusion for $g \in C^{\infty}(\partial C_{\eps})$ since  the constant $C$ in the variant of \eqref{fund-est-L2} for $C_{\eps, \delta}$ does not depend on $\eps$ and $\delta$. The general result stated in Lemma~\ref{lem-est-L2} follows by the standard approximate arguments.    

\medskip
We now give the proof for the case $v \in C^{\infty}(\overline{B_{4} \setminus C_{\eps}})$. Applying Lemma~\ref{lem-Morawetz} and using the fact that $\langle x, \eta \rangle>  0$ on $\partial C_{\eps}$, we obtain 
\begin{multline}\label{toto-1}
\frac{1}{2} \int_{B_{4} \setminus C_{\eps}} \big( |\nabla v|^{2} + k^{2} |v|^{2} \big)  +  \int_{\partial C_{\eps}} \Re \Big(  \langle \nabla v , \eta \rangle 
\langle \nabla \bar v , x \rangle \Big) 
- \frac{1}{2}  \int_{\partial C_{\eps}}   \langle x, \eta \rangle |\nabla v|^{2} \\[6pt]
\le  \int_{\partial C_{\eps}} \Big| \frac{\partial v}{\partial \eta} \bar v \Big|   + C \Big(  \| \nabla v\|_{L^{2}(\partial B_{4})}^{2} +  \| v\|_{L^{2}(\partial B_{4})}^{2}\Big).
\end{multline}
We recall that $\eta$ denotes the exterior normal unit vector on $\partial C_{\eps}$; this is the reason for which one has the change in sign for the second and third integrals in comparison with Lemma~\ref{lem-Morawetz}. We have, on $\partial C_{\eps}$,
\begin{align}\label{toto0}
\Re \Big(  \langle \nabla v , \eta \rangle 
\langle \nabla \bar v , x \rangle \Big) 
- \frac{1}{2}  \langle x, \eta \rangle |\nabla v|^{2} & \ge  \frac{1}{2} \langle x, \eta \rangle |\partial_{\eta} v|^{2} - \frac{1}{2} \langle x, \eta \rangle |\nabla_{\Gamma} v|^{2} -  |\langle \nabla v, \alpha (x) \rangle| \cdot |\partial_{\eta} v|  \nonumber\\[6pt]
& \ge  \frac{1}{4} \langle x, \eta \rangle |\partial_{\eta} v|^{2} - \frac{1}{2} \langle x, \eta \rangle |\nabla_{\Gamma} v|^{2} -  \frac{1}{\langle x , \eta \rangle} |\langle \nabla v, \alpha (x) \rangle|^{2}, 
\end{align}
where $\alpha(x): = x - \langle x , \eta \rangle \eta$. Here in the last inequality, we used the Cauchy inequality:  
\begin{equation*}
 |\langle \nabla v, \alpha (x) \rangle| \cdot |\partial_{\eta} v|  \le \frac{1}{\langle x, \eta \rangle}|\langle \nabla v, \alpha (x) \rangle|^{2} + \frac{1}{4} \langle x, \eta \rangle |\partial_{\eta} v|^{2}. 
\end{equation*}
We also have, on $\partial C_{\eps}$, 
\begin{equation}\label{toto1}
\Big| \frac{\partial v}{\partial \eta} \bar v \Big| \le \frac{1}{4} \langle x, \eta \rangle |\partial_{\eta} v|^{2} + \frac{1}{\langle x, \eta \rangle} |v|^{2}. 
\end{equation}
It follows from \eqref{toto-1}, \eqref{toto0}, and \eqref{toto1} that 
\begin{multline*}
 C \left\{\int_{B_{4} \setminus C_{\eps}} \big( |\nabla v|^{2} + k^{2} |v|^{2} \big) +  \int_{\partial C_{\eps}}  \langle x, \eta \rangle |\partial_{\eta}  v|^{2} \right\}
\le   \| \nabla v\|_{L^{2}(\partial B_{4})}^{2} +  \| v\|_{L^{2}(\partial B_{4})}^{2} \\[6pt]
+ \int_{\partial C_{\eps}}  \Big( \langle x, \eta \rangle |\nabla_{T} v|^{2} + |\langle \nabla v, \alpha (x) \rangle| \cdot |\partial_{\eta} v|  \Big)+ \int_{\partial C_{\eps}}  \frac{1}{\langle x, \eta \rangle}  |v|^{2}. 
\end{multline*}
Since $\langle \alpha(x), \eta \rangle = 0$, it follows that
\begin{equation*}
 |\langle \nabla v, \alpha (x) \rangle| \cdot |\partial_{\eta} v|  \le \frac{8|x|^{2}}{C \langle x, \eta \rangle } |\nabla_{\Gamma} v|^{2}
 + \frac{C}{2} \langle x, \eta \rangle |\partial_{\eta } v|^{2}. 
\end{equation*}
Therefore, the conclusion follows. 
\proofend

\medskip
The following lemma will be used in the proof of Theorem~\ref{thm1}. 

\begin{lemma}\label{lem-stability-k} Let $0 < \eps < 1$,  $f_\eps \in H^{-1}(\mR^3 \setminus C_{\eps})$  with $\supp f_\eps \subset B_3 \setminus C_{\eps}$ and let $u_\eps \in H^1_{\loc}(\mR^3 \setminus C_\eps)$ be the unique outgoing solution to the equation
\begin{equation*}
\Delta u_\eps + k^2 u_\eps = f_\eps \mbox{ in } \mR^3 \setminus C_\eps \quad \mbox{ with } \quad u_\eps =0 \mbox{ on } \partial C_\eps.
\end{equation*}
Then
\begin{equation}\label{stability}
\| u_\eps \|_{H^1(B_5 \setminus C_\eps)} \le C \|f_\eps \|_{H^{-1}},
\end{equation}
\footnote{$\| f_{\eps}\|_{H^{-1}}: = \mathop{\sup}_{ \varphi \in C^{1}_{c}(\mR^{3} \setminus C_{\eps}) ; \; \| \varphi \|_{H^{1}(\mR^{3} \setminus C_{\eps})} =1} |\langle f_{\eps}, \varphi \rangle |$} for some positive constant $C$ depending on $k$ but independent of $f_\eps$ and $\eps$.
\end{lemma}

\noindent{\bf Proof.} We first prove
\begin{equation}\label{toto-claim}
\| u_\eps \|_{L^2(B_5 \setminus C_\eps)} \le C \|f_\eps \|_{H^{-1}},
\end{equation}
for some positive constant $C$ independent of $f_\eps$ and $\eps$ by contradiction. Suppose that this is not true. Then there exist a sequence $\eps_n \to 0_+$, a sequence $(f_n) \subset H^{-1}(\mR^3 \setminus C_{\eps_{n}})$ with $\supp f_n \subset B_3 \setminus C_{\eps_{n}}$, and a sequence $(u_n) \subset H^1_{\loc}(\mR^3 \setminus C_{\eps_n})$, the unique
outgoing solution to the equation
\begin{equation*}
\Delta u_n + k^2 u_n = f_n \mbox{ in } \mR^3 \setminus C_{\eps_n} \quad \mbox{ and } \quad u_n = 0 \mbox{ on } \partial C_{\eps_n},
\end{equation*}
such that
\begin{equation}\label{contra}
\| u_n \|_{L^2(B_5 \setminus C_{\eps_{n}})} = 1 \quad \mbox{ and } \| f_n \|_{H^{-1}}  \to 0.
\end{equation}
Multiplying the equation of $u_n$ by $\bar u_n$ (the conjugate of $u_n$) and integrating on  $B_4 \setminus C_{\eps_{n}}$, we have
\begin{equation}\label{mu1}
- \int_{B_4 \setminus C_{\eps_{n}}} |\nabla u_n|^2 + k^2 \int_{B_4 \setminus C_{\eps_{n}}} |u_n|^2 + \int_{\partial B_4} \partial_r u_n \bar u_n = \int_{B_4 \setminus C_{\eps_{n}}} f_n \bar u_n.
\end{equation}
Since $\Delta u_{n} + k^{2} u_{n} = 0$ in $B_{5} \setminus B_{3}$ and $\| u_{n}\|_{L^{2}(B_{5} \setminus B_{3})} \le 1$, it follows from the standard estimates for the Laplace equation that $\| u_{n}\|_{C^{1}(B_{9/2} \setminus B_{7/2})} \le C$; which yields
\begin{equation}\label{mu2}
\Big| \int_{\partial B_4} \partial_r u_n \bar u_n \Big| \le C.  
\end{equation}
Since 
$$
\left| \int_{B_4 \setminus C_{\eps_{n}}} f_n \bar u_n \right| \le C \| f_{n}\|_{H^{-1}} \| u_{n}\|_{H^{1}(B_{4} \setminus C_{\eps_{n}})}, 
$$
we derive from \eqref{contra}, \eqref{mu1} and \eqref{mu2} that 
\begin{equation}\label{e-der}
\int_{B_4 \setminus C_{\eps_{n}}} |\nabla u_n|^2 \le C.
\end{equation}
By the representation formula for the Helmholtz equation in $3d$, it follows from \eqref{contra} and \eqref{e-der} that
\begin{equation}\label{toto}
\sup_n \|u_n\|_{H^1(B_r \setminus C_{\eps_{n}})}  \le C(r) \quad \forall \, r > 0. 
\end{equation}
Extend $u_n$ by 0 in $C_{\eps_n}$ and still denote this extension by $u_{n}$. We have
\begin{equation*}
\sup_{n} \| u_n\|_{H^1(B_r)} < C(r) \quad \forall \, r > 0.
\end{equation*}
Without loss of generality, one might assume that $u_{n} \to u$ in $H^{1}_{\loc}(\mR^{3})$. Moreover, $u$ satisfies the outgoing condition and 
\begin{equation}\label{eq-lim}
\Delta u + k^2 u = 0 \mbox{ in } \mR^3 \setminus \{(0, 0, z); |z| \le 1/2 \}.
\end{equation}

We claim that $u$ satisfies the Helmholtz equation in the whole space, i.e.,  
\begin{equation}\label{claim}
\Delta u + k^2 u = 0 \mbox{ in } \mR^3.
\end{equation}

Indeed, for $\delta \ge 0$, set 
\begin{equation*}
\Omega_{\delta} = \{x \in \mR^{3}; \dist (x, (0, 0, z)) \le \delta \mbox{ for } |z| \le 1/2\}, 
\end{equation*}
and for $\delta > 0$, let $\phi \in C^1(\mR^3)$ be such that
$\phi(x)  =1 $ if $\dist(x, \Omega_{0}) \le \delta$, $\phi(x) = 0$ if $\dist(x, \Omega_{0}) \ge 2 \delta$, $0 \le \phi \le 1$,  and $|\nabla \phi| \le C/ \delta$ for some fixed positive constant $C$. For $\varphi \in C^1_{\mc}(\mR^3)$, define
\begin{equation*}
 \varphi_1 = \varphi (1 - \phi)
\mbox{ and } \varphi_2 = \varphi \phi 
\end{equation*}
We have
\begin{equation*}
\int_{\mR^3} \nabla u \nabla \varphi - k^2 \int_{\mR^3} u \varphi = \int_{\mR^3} \nabla u \nabla \varphi_1 - k^2 \int_{\mR^3} u \varphi_1 +
\int_{\mR^3} \nabla u \nabla \varphi_2 - k^2 \int_{\mR^3} u \varphi_2.
\end{equation*}
Since $\supp \varphi_{1} \subset \mR^{3} \setminus \Omega_{\delta}$, it follows from  \eqref{eq-lim} and the condition on the support of $\phi$ that 
\begin{equation}\label{mu3}
\int_{\mR^3} \nabla u \nabla \varphi - k^2 \int_{\mR^3} u \varphi = \int_{\Omega_{2 \delta}} \nabla u \nabla \varphi_2 - k^2 \int_{\Omega_{2 \delta}} u \varphi_2.
\end{equation}
On the other hand,
\begin{equation}\label{mu4}
\Big| \int_{\Omega_{2 \delta}} \nabla u \nabla \varphi_2 - k^2 \int_{\Omega_{2 \delta}} u \varphi_2 \Big| \le C \int_{\Omega_{2\delta}} \delta^{-1}(|\nabla u| + |u|). 
\end{equation}
By H\"older's inequality, 
\begin{equation}\label{mu5}
\int_{\Omega_{2\delta}}(|\nabla u| + |u|) \le  C\delta \Big(\int_{\Omega_{2\delta}}(|\nabla u|^{2} + |u|^{2}) \Big)^{1/2} .
\end{equation}
A combination of \eqref{mu3}, \eqref{mu4}, and \eqref{mu5} yields 
\begin{equation*}
\Big| \int_{\mR^3} \nabla u \nabla \varphi - k^2 \int_{\mR^3} u \varphi  \Big| \le C \Big(\int_{\Omega_{2\delta}}(|\nabla u|^{2} + |u|^{2}) \Big)^{1/2} \to 0 \mbox{ as } \delta \to 0. 
\end{equation*}
Hence \eqref{claim} is proved. 

\medskip

We derive from \eqref{claim} that  $u=0$ since $u$ satisfies the outgoing condition. We have a contradiction. Therefore, 
\begin{equation*}
\| u_{\eps} \|_{L^{2}(B_{5} \setminus B_{\eps})} \le C \| f\|_{H^{-1}}:
\end{equation*} 
\eqref{toto-claim} is proved. 
Multiplying the equation of $u_\eps$ by $\bar u_\eps$ and integrating on $B_4 \setminus C_\eps$, we have
\begin{equation*}
- \int_{B_4 \setminus C_{\eps}} |\nabla u_\eps|^2 + k^2 \int_{B_4 \setminus C_{\eps}} |u_\eps|^2 + \int_{\partial B_4} \partial_r u_\eps \bar u_\eps = \int_{B_4 \setminus C_{\eps}} f_\eps \bar u_\eps.
\end{equation*}
Similar to \eqref{toto}, we obtain
\begin{equation*}
\|u_\eps\|_{H^1(B_r \setminus C_\eps)}  \le C(r) \| f_{\eps}\|_{H^{-1}}.
\end{equation*}
The proof is complete.
\proofend

\medskip

As a consequence of Lemma~\ref{lem-stability-k}, we have
\begin{corollary}\label{cor1}
Let $k>0$, $0 < \eps < 1$,  $g_\eps \in C^{2}(C_{3 \eps} \setminus C_{\eps})$, and let $u_\eps \in H^1_{\loc}(\mR^3 \setminus C_\eps)$ be the unique outgoing solution to the equation
\begin{equation*}
\Delta u_\eps + k^2 u_\eps = 0 \mbox{ in } \mR^3 \setminus C_\eps \quad \mbox{ with } \quad u_\eps = g_{\eps} \mbox{ on } \partial C_\eps.
\end{equation*}
Then
\begin{equation}\label{stability}
\| u_\eps \|_{H^1(B_5 \setminus C_\eps)} \le C \big( \| g_{\eps}\|_{C^{0}} + \eps \|\nabla g_{\eps}\|_{C^{0}} \big),
\end{equation}
for some positive constant $C$ depending on $k$ but independent of $g_\eps$ and $\eps$.
\end{corollary}

\noindent{\bf Proof.} Let $\varphi \in C^{1}(\mR^{3})$  be such that $\varphi(x) = 1$ if $\dist(x, C_{\eps}) \le \eps$, $\varphi(x) = 0$ if $\dist(x, C_{\eps}) \ge 2 \eps$, $|\varphi| \le 1$ and $|\nabla \varphi| \le C /\eps$, for some fixed constant $C$. 
Set 
\begin{equation*}
v_{\eps} = u_{\eps} -  \varphi g_{\eps}.  
\end{equation*}
We have 
\begin{equation*}
\Delta v_{\eps} + k^{2} v_{\eps} = f_{\eps} \mbox{ in } \mR^{3} \setminus C_{\eps} \quad \mbox{ and } \quad v_{\eps} = 0 \mbox{ on } \partial C_{\eps},
\end{equation*}
where 
\begin{equation*}
f_{\eps} = -\Delta (\varphi g_{\eps}) - k^{2} \varphi g_{\eps}. 
\end{equation*}
We have 
\begin{equation*}
\| f_{\eps}\|_{H^{-1}} \le C\Big( \| \nabla (\varphi g_{\eps}) \|_{L^{2}} + k^{2} \| \varphi g_{\eps}\|_{L^{2}} \Big) \le C   \big( \| g_{\eps}\|_{C^{0}} + \eps \| \nabla g_{\eps}\|_{C^{0}} \big). 
\end{equation*}
The conclusion now follows from Lemma~\ref{lem-stability-k}. \proofend

\subsection{Proof of Proposition~\ref{pro1}}

Set
\begin{equation*}
v_\eps = u_\eps - u.
\end{equation*}
Then $v_\eps \in H^1_{\loc}(\mR^3 \setminus C_\eps)$ is the unique outgoing solution to the equation
\begin{equation*}
\Delta v_\eps + k^2 v_\eps = 0 \mbox{ in } \mR^3 \setminus \bar C_\eps \quad \mbox{ with } \quad v_\eps = - u \mbox{ on } \partial C_\eps.
\end{equation*}
Let $w_{1, \eps} \in H^1_{\loc}(\mR^3 \setminus C_\eps)$ be the unique outgoing solution to the equation
\begin{equation}\label{w1}
\Delta w_{1,\eps} + k^2 w_{1, \eps} = 0 \mbox{ in } \mR^3 \setminus \bar C_\eps \quad \mbox{ with } \quad w_{1, \eps} (x', z) = - u(0, z) \mbox{ on } \partial C_\eps.
\end{equation}
and $w_{2, \eps} \in H^1_{\loc}(\mR^3 \setminus C_\eps)$ be the unique outgoing solution to the equation
\begin{equation}\label{w2}
\Delta w_{2, \eps} + k^2 w_{2, \eps} = 0 \mbox{ in } \mR^3 \setminus \bar C_\eps \quad \mbox{ with } \quad w_{2, \eps}(x', z) = - u(x', z) + u(0, z) \mbox{ on } \partial C_\eps.
\end{equation}
It is clear that
\begin{equation*}
v_\eps = w_{1, \eps} + w_{2, \eps}.
\end{equation*}
Set
\begin{equation*}
\Gamma_{\eps, s} = \partial C_r \cap \{z =  - 1/2, 1/2\},
\end{equation*}
and
\begin{equation*}
\Gamma_{\eps, c} = \partial C_r \cap \{ |z| <1/2\}.
\end{equation*}
We have, from Lemma~\ref{lem-est-L2}, 
\begin{equation}\label{est1-k}
\int_{\Gamma_{\eps, s}} |\partial_\eta w_{1, \eps}|^2 + \eps \int_{\Gamma_{\eps, c}} |\partial_\eta w_{1, \eps}|^2  \le C  \|f \|_{L^2}^{2}, \footnote{Roughly speaking, $\partial_{\eta} w_{1, \eps}(x)$ is of  order $1/\eps$ on $\partial C_{\eps}$ as mentioned in the Introduction.}
\end{equation}
and, from Corollary~\ref{cor1},  
\begin{equation}\label{est2-k}
\| w_{2, \eps}\|_{H^{1}(B_{5} \setminus B_{2})} \le C \eps \|f \|_{L^2}.
\end{equation}
Here we also applied the standard regularity theory of  elliptic equations to obtain 
\begin{equation*}
\| u\|_{C^{1}(\bar B_{1})} \le C \| f\|_{L^{2}},
\end{equation*}
since $\Delta u + k^{2}u = 0$ in $B_{3/2}$ and $\| u \|_{L^{2}(B_{3/2})} \le C \| f\|_{L^{2}}$.  

By the representation formula (see e.g. \cite{ColtonKressInverse}), we have
\begin{equation*}
w_{1, \eps}(x) =  \int_{\partial C_\eps} G_k(x, y) \frac{\partial w_{1, \eps}}{\partial \eta}(y)  \, d \Gamma(y) -  \int_{\partial C_\eps} \frac{\partial G_k(x,y)}{\partial \eta } w_{1, \eps} (y) \, d \Gamma(y). 
\end{equation*}
Here $G_{k}(x,y)$ is the fundamental solution to the Helmholtz equation, i.e.,  $G_{k}(x, y) = \frac{e^{ ik|x-y|}}{4 \pi |x - y|}$. 
It follows that 
\begin{multline}\label{est1.1-k}
w_{1, \eps}(x) =  \int_{\partial C_\eps} [G_k(x, y) - G_k(\pi(x), y)] \frac{\partial w_{1, \eps}}{\partial \eta}(y)  \, d \Gamma(y) -  \int_{\partial C_\eps} \frac{\partial G_k(x,y)}{\partial \eta } w_{1, \eps} (y) \, d \Gamma(y) \\[6pt]
+ \int_{\partial C_\eps}  G_k(\pi(x), y) \frac{\partial w_{1, \eps}}{\partial \eta}(y)  \, d \Gamma(y),
\end{multline}
where $\pi(x) = (0, 0, x_3)$ for $x=(x_1, x_2, x_3)$. Applying Lemma~\ref{lem-sym} for $w_{1, \eps}$, we have
\begin{equation}\label{est1.2-k}
\int_{\partial C_\eps}  G_k(\pi(x), y) \frac{\partial w_{1, \eps}}{\partial \eta}(y)  \, d \Gamma(y) = \int_{\Gamma_{\eps, s}}  G_k(\pi(x), y) \frac{\partial w_{1, \eps}}{\partial \eta}(y)  \, d \Gamma(y), 
\end{equation}
since $G_k(\pi(x), y)$ is constant on $\partial C_{\eps} \cap \{(x', z); \; z = a \}$ for $a \in (-1/2, 1/2)$. 
A combination of  \eqref{est1-k}, \eqref{est1.1-k} and  \eqref{est1.2-k} yields
\begin{equation}\label{est3-k}
|w_{1, \eps}| \le C \eps \|f \|_{L^2} \quad \mbox{ for } x \in B_{5} \setminus B_{2}.
\end{equation}
Since $v_{\eps} = w_{1, \eps} + w_{2, \eps}$, the conclusion follows from \eqref{est2-k} and \eqref{est3-k}. \proofend

\section{Proof of Proposition~\ref{pro2}}\label{sect-Pro2}

The following result (see e.g., \cite{NguyenHelmholtz, NguyenVogelius2, NguyenPerfectCloaking}) will be used in the proof of Proposition~\ref{pro2}. 

\begin{lemma}\label{lemNg1} Let $k, \tau > 0$,  $0 < \eps < k$, $D \subset B_1$ be a smooth open connected subset of $\mR^d$ $d=2, \, 3$, and $g \in H^{\frac{1}{2}}(\partial D)$. Assume that $\mR^d \setminus D$ is connected and $v_\eps \in H^1_{\loc}(\mR^3 \setminus D)$ is the unique outgoing solution to the system
\begin{equation*}
\left\{\begin{array}{ll}
\Delta v_\eps + \eps^2  v_\eps = 0 & \mbox{in } \mR^d \setminus D, \\[6pt]
v_\eps = g & \mbox{on } \partial D. 
\end{array}\right.
\end{equation*}
Then
\begin{equation*}
\| v_\eps \|_{H^1(B_5 \setminus D)} \le C \| g\|_{H^\frac{1}{2}(\partial D)}
\end{equation*}
and, for $x$ with $|x| \in (\tau/\eps, 6 \tau/ \eps)$, 
\begin{equation*}
|v_{\eps}(x)| \le \left\{ \begin{array}{cl}
\dsp C(\tau) \eps \| g\|_{H^\frac{1}{2}(\partial D)} & \mbox{ if } d= 3, \\[6pt]
\dsp \frac{C(\tau)}{|\ln \eps |}\| g\|_{H^\frac{1}{2}(\partial D)} & \mbox{ if } d =2.  
\end{array} \right.
\end{equation*}
for some positive constants $C$ and $C(\tau)$ independent of $g$ and $\eps$. 
\end{lemma}

We are ready to give

\medskip
\noindent {\bf Proof of Proposition~\ref{pro2}.} Define
\begin{equation*}
w_{\eps} = u_{\eps} - u. 
\end{equation*}
It follows that $w_{\eps} \in H^1_{\loc}(\mR^d)$ is the unique outgoing solution to the system
\begin{equation}
\left\{\begin{array}{ll}
\Delta w_{ \eps} + k^2 w_{\eps} = 0 & \mbox{ in } \mR^d \setminus B_\eps,\\[6pt]
w_{ \eps} = - u & \mbox{ in } \partial B_\eps. 
\end{array}\right.
\end{equation}
Define $W_{\eps}(x) = w_{\eps}(\eps x)$. Then $W_{\eps} \in H^1_{\loc}(\mR^d)$ is the unique outgoing solution to the system
\begin{equation}
\left\{\begin{array}{ll}
\Delta W_{\eps} + \eps^2 k^2 W_{\eps} = 0 & \mbox{ in } \mR^d \setminus B_{1},\\[6pt]
W_{\eps} = -  u (\eps \cdot) & \mbox{ in } \partial B_{1}. 
\end{array}\right.
\end{equation}
Let  $W_{1,  \eps}$ be the unique outgoing solution to the equation 
\begin{equation}
\left\{\begin{array}{ll}
\Delta W_{1,  \eps} + \eps^2 k^2 W_{1, \eps} = 0 & \mbox{ in } \mR^d \setminus B_{1},\\[6pt]
W_{1,  \eps} = -  u (0) & \mbox{ in } \partial B_{1}, 
\end{array}\right.
\end{equation}
and $W_{2, \eps} = W_{\eps} - W_{1, \eps}$. Then  $W_{2, \eps}$ is the unique outgoing solution to the system
\begin{equation}
\left\{\begin{array}{ll}
\Delta W_{2,  \eps} + \eps^2 k^2 W_{2,  \eps} = 0 & \mbox{ in } \mR^d \setminus B_{1},\\[6pt]
W_{2, \eps} = -  u (\eps \cdot) + u(0) & \mbox{ in } \partial B_{1}. 
\end{array}\right.
\end{equation}
We will only consider the $3d$ case; the $2d$ case follows similarly. By Lemma~\ref{lemNg1}, it follows that, for $1/ \eps \le |x| \le 6/ \eps$,
\begin{equation}\label{est1}
|W_{2, \eps}| \le C \eps^{2} \|f \|_{L^{2}}. 
\end{equation}
Here we applied the standard regularity theory  of elliptic equations to obtain 
\begin{equation*}
\|u(\eps \cdot) - u \|_{C^{1}(B_{2})} \le C \eps \| f\|_{L^{2}}. 
\end{equation*}

We now estimate $W_{1, \eps}$. We have, by the representation formula,  
\begin{equation*}
W_{1, \eps} (x) = \int_{\partial B_{1}} G_{\eps k}(x, y) \frac{\partial W_{1, \eps}}{\partial r} (y) \, d \Gamma (y) - \int_{\partial B_{1}} \frac{ \partial G_{\eps k}(x,y)}{\partial r}  W_{1, \eps} (y) \, d \Gamma (y). 
\end{equation*}
Here $G_{\eps k}(x,y)$ is the fundamental solution to the Helmholtz equation i.e. $G_{\eps k}(x, y) = \frac{e^{ i\eps k|x-y|}}{4 \pi |x - y|}$. 
Since 
\begin{equation*}
\int_{\partial B_{1}} \frac{\partial W_{1, \eps}}{\partial r}(y) d \Gamma(y) = 0, 
\end{equation*}
it follows that 
\begin{equation*}
W_{1, \eps} (x) = \int_{\partial B_{1}} \big[G_{\eps k}(x, y) - G_{\eps k }(x, 0)\big]\frac{\partial W_{1, \eps}}{\partial r} (y) \, dy - \int_{\partial B_{1}} \frac{ \partial G_{\eps k} (x, y)}{\partial r}  W_{1, \eps} (y) \, dy. 
\end{equation*}
For $1/\eps \le  |x| \le 6/ \eps$, we have
\begin{equation*}
\sup_{y \in \partial B_{1}}\left|G_{\eps k}(x, y) - G_{\eps k}(x, 0) \right| + \sup_{y \in \partial B_{1}} \left| \frac{\partial G_{\eps k}(x, y)}{\partial r}  \right| \le C \eps^{2}. 
\end{equation*}
On the other hand, 
\begin{equation*}
\sup_{y \in \partial B_{1}}\left| \frac{\partial W_{1, \eps}}{\partial r} (y) \right|  + \sup_{y \in \partial B_{1}}\left| W_{1, \eps} (y) \right| \le C | u(0) | \le C \|f \|_{L^{2}}. 
\end{equation*}
In the last inequality, we used the standard regularity theory of elliptic equations.
It follows that, for $1/ \eps \le |x| \le 6/ \eps$,  
\begin{equation}\label{est2}
|W_{1, \eps}(x)| \le C \eps^{2} \| f\|_{L^{2}}. 
\end{equation}
A combination of \eqref{est1} and \eqref{est2} implies 
\begin{equation*}
\| w_{\eps}\|_{L^{2}(B_{6} \setminus B_{1})} \le C \eps^{2} \|f \|. 
\end{equation*}
Since $\Delta w_{\eps} + k^{2} w_{\eps} = 0$ in $B_{6} \setminus B_{1}$, it follows from standard estimates for the Laplace equation that 
\begin{equation*}
\| w_{\eps}\|_{H^{1}(B_{5} \setminus B_{2})} \le C \eps^{2} \|f \|. 
\end{equation*}
The proof is complete. 
\proofend




\providecommand{\bysame}{\leavevmode\hbox to3em{\hrulefill}\thinspace}
\providecommand{\MR}{\relax\ifhmode\unskip\space\fi MR }
\providecommand{\MRhref}[2]{%
  \href{http://www.ams.org/mathscinet-getitem?mr=#1}{#2}
}
\providecommand{\href}[2]{#2}


\begin{thebibliography}{10}

\bibitem{AluEngheta}
A.~Alu and N.~Engheta, \emph{{Achieving transparency with plasmonic and
  metamaterial coatings}}, Phys. Rev. E \textbf{95} (2005), 106623.

\bibitem{AluEngheta2}
\bysame, \emph{{Multifrequency Optical Invisibility Cloak with Layered
  Plasmonic Shells}}, Phys. Rev. Lett. \textbf{100} (2008), 113901.

\bibitem{AmmariCiraoloKangLeeMilton}
H.~Ammari, G.~Ciraolo, H.~Kang, H.~Lee, and G.~W. Milton, \emph{{Spectral
  theory of a Neumann-Poincar\'e-type operator and analysis of cloaking due to
  anomalous localized resonance}},  (2011), preprint.

\bibitem{AmmariKang-Po}
H.~Ammari and H.~Kang, \emph{{Polarization and Moment Tensors: With
  Applications to Inverse Problems and Effective Medium Theory}}, Applied
  Mathematical Sciences, vol. 162, Springer, New York, 2007.

\bibitem{AmmariKHL-Helmholtz}
H.~Ammari, H.~Kang, H.~Lee, and M.~Lim, \emph{{Enhancement of near-cloaking.
  Part II: the Helmholtz equation}},  (2011), preprint.

\bibitem{AmmariKHL-Conductivity}
\bysame, \emph{{Enhancement of Near Cloaking Using Generalized Polarization
  Tensors Vanishing Structures. Part I: The Conductivity Problem}},  (2011),
  preprint.

\bibitem{BouchitteSchweizer10}
G.~Bouchitt\'e and B.~Schweizer, \emph{{Cloaking of small objects by anomalous
  localized resonance}}, Quart. J. Mech. Appl. Math. \textbf{63} (2010),
  437--463.

\bibitem{Kildishev07}
W.~Cai, U.~K. Chettiar, A.~V. Kildishev, and V.~M. Shalaev, \emph{{Optical
  cloaking with metamaterials}}, Nature Photonics \textbf{1} (2007), 224--227.

\bibitem{Capdeboscq}
Y.~Capdeboscq, \emph{{On the scattered field generated by a ball imhomogeneity
  of constant index}},  (2012), preprint.

\bibitem{CapdeboscqVogelius}
Y.~Capdeboscq and M.~Vogelius, \emph{{ A general representation formula for
  boundary voltage perturbations caused by internal conductivity
  inhomogeneities of low volume fraction}}, M2AN Math. Model. Numer. Anal.
  \textbf{37} (2003), 159--173.

\bibitem{Cedio-FengyaMoskowVogelius}
D.~J. Cedio-Fengya, S.~Moskow, and M.~S. Vogelius, \emph{{Identification of
  conductivity imperfections of small diameter by boundary measurements.
  Continuous dependence and computational reconstruction}}, Inverse Problems
  \textbf{14} (1999), 553.

\bibitem{ColtonKressInverse}
D.~Colton and R.~Kress, \emph{{Inverse acoustic and electromagnetic scattering
  theory}}, second ed., Applied Mathematical Sciences, vol.~98,
  Springer-Verlag, Berlin, 1998.

\bibitem{FriedmanVogelius}
A.~Friedman and M.~Vogelius, \emph{{Identification of small inhomogeneities of
  extreme conductivity by boundary measurements: a theorem on continuous
  dependence}}, Arch. Rational Mech. Anal. \textbf{105} (1989), 299--326.

\bibitem{GreenleafKurylevLassasLeonhardtUhlmann}
A.~Greenleaf, Y.~Kurylev, M.~Lassas, U.~Leonhard, and G.~Uhlmann,
  \emph{{Cloaked electromagnetic, acoustic, and quantum amplifiers via
  transformation optics}}, Proc. Natl. Acad. Sci. \textbf{109} (2012),
  10169--10174.

\bibitem{GreenleafKurylevLassasUhlmann-QP}
A.~Greenleaf, Y.~Kurylev, M.~Lassas, and G.~Uhlmann, \emph{{Approximate quantum
  cloaking and almost trapped states}}, Phys. Rev. Lett.

\bibitem{GreenleafKurylevLassasUhlmann07}
\bysame, \emph{{Full-wave invisibility of active devices at all frequencies}},
  Comm. Math. Phys. \textbf{275} (2007), 749--789.

\bibitem{GreenleafKurylevLassasUhlmann07SHS}
\bysame, \emph{{Improvement of cylindrical cloaking with the SHS lining}}, Opt.
  Exp. \textbf{15} (2007), 12717.

\bibitem{GreenleafKurylevLassasUhlmann-Q}
\bysame, \emph{{Approximate Acoustic and Quantum Cloaking}}, Jour. Spectral
  Theory \textbf{1} (2011), 27--80.

\bibitem{GreenleafKurylevLassasUhlmann-PE}
\bysame, \emph{{Cloaking a Sensor via Transformation Optics}}, Phys. Rev. E
  \textbf{83} (2011), 016603.

\bibitem{GreenleafLassasUhlmann}
A.~Greenleaf, M.~Lassas, and G.~Uhlmann, \emph{{On nonuniqueness for Calderon's
  inverse problem}}, Math. Res. Lett. \textbf{10} (2003), 685--693.

\bibitem{HaddarJolyNguyen1}
H.~Haddar, P.~Joly, and H-M. Nguyen, \emph{{Generalized impedance boundary
  conditions for scattering by strongly absorbing obstacles: the scalar case}},
  Math. Models Methods Appl. Sci. \textbf{15} (2005), 1273--1300.

\bibitem{KohnOnofreiVogeliusWeinstein}
R.~V. Kohn, D.~Onofrei, M.~S. Vogelius, and M.~I. Weinstein, \emph{{Cloaking
  via change of variables for the Helmholtz equation}}, Comm. Pure Appl. Math.
  \textbf{63} (2010), 973--1016.

\bibitem{KohnShenVogeliusWeinstein}
R.~V. Kohn, H.~Shen, M.~S. Vogelius, and M.~I. Weinstein, \emph{{Cloaking via
  change of variables in electric impedance tomography}}, Inverse Problems
  \textbf{24} (2008), 015016.

\bibitem{LaiChenZhangChanComplementary}
Y.~Lai, H.~Chen, Z.~Zhang, and C.~T. Chan, \emph{{Complementary Media
  Invisibility Cloak that Cloaks Objects at a Distance Outside the Cloaking
  Shell}}, Phys. Rev. Lett. \textbf{102} (2009).

\bibitem{LassasZhou}
M.~Lassas and T.~Zhou, \emph{{Two dimensional invisibility cloaking for
  Helmholtz equation and non-local boundary conditions}}, Math. Res. Lett.
  \textbf{18} (2011), 473--488.

\bibitem{Leonhardt}
U.~Leonhardt, \emph{{Optical conformal mapping}}, Science \textbf{312} (2006),
  1777--1780.

\bibitem{LeonhardtTyc}
U.~Leonhardt and T.~Tyc, \emph{{Broadband Invisibility by Non-Euclidean
  Cloaking}}, Science \textbf{323} (2009), 110--112.

\bibitem{Liu}
H.~Liu, \emph{{Virtual reshaping and invisibility in obstacle scattering}},
  Inverse Problems \textbf{25} (2009), 045006.

\bibitem{LiuSun}
H.~Liu and H.~Sun, \emph{{Enhanced near-cloak by FSH lining}},  (2012),
  preprint.

\bibitem{LiuZhou}
H.~Y. Liu and T.~Zhou, \emph{{Two Dimensional Invisibility Cloaking via
  Transformation Optics}}, Discrete Contin. Dyn. Syst. \textbf{31} (2011),
  525--543.

\bibitem{Miller}
D.~Miller, \emph{{On perfect cloaking}}, Optics Express \textbf{14} (2006),
  12457.

\bibitem{MiltonNicorovici}
G.~W. Milton and N-A.~P. Nicorovici, \emph{{On the cloaking effects associated
  with anomalous localized resonance}}, Proc. R. Soc. Lond. Ser. A \textbf{462}
  (2006), 3027--3059.

\bibitem{Milton-folded}
G.~W. Milton, N.~P. Nicorovici, R.~C. McPhedran, K.~Cherednichenko, and
  Z.~Jacob, \emph{{Solutions in folded geometries, and associated cloaking due
  to anomalous resonance}}, New J. Phys. \textbf{10} (2008), 115021.

\bibitem{MiltonNicoroviciMcPhedranPodolskiy}
G.~W. Milton, N.~P. Nicorovici, R.~C. McPhedran, and V.~A. Podolskiy, \emph{{A
  proof of superlensing in the quasistatic regime and limitations of
  superlenses in this regime due to anomalous localized resonance }}, Proc. R.
  Soc. Lond. Ser. A \textbf{461} (2005), 3999--4034.

\bibitem{MorawetzLudwig}
C.~S. Morawetz and D.~Ludwig, \emph{{An inequality for the reduced wave
  operator and the justification of geometrical optics}}, Comm. Pure Appl.
  Math. \textbf{21} (1968), 187--203.

\bibitem{NguyenHelmholtz}
H-M. Nguyen, \emph{{Cloaking for the Helmholtz equation in the whole spaces}},
  Comm. Pure Appl. Math. \textbf{63} (2010), 1505--1524.

\bibitem{NguyenNegative}
\bysame, \emph{{A study of negative index materials using transformation optics
  with applications to super lenses, cloaking, and illusion optics: the scalar
  case}},  (2012), preprint.

\bibitem{NguyenPerfectCloaking}
\bysame, \emph{{Approximate cloaking for the Helmholtz equation via
  transformation optics and consequences for perfect cloaking}}, Comm. Pure
  Appl. Math. \textbf{65} (2012), 155--186.

\bibitem{NguyenVogelius}
H-M. Nguyen and M.~S. Vogelius, \emph{{A representation formula for the voltage
  perturbations caused by diametrically small conductivity inhomogeneities.
  Proof of uniform validity}}, Ann. Inst. H. Poincar\'e Anal. Non Lin\'eaire
  \textbf{26} (2009), 2283--2315.

\bibitem{NguyenVogelius3}
\bysame, \emph{{Approximate cloaking for the wave equation via change of
  variables}}, SIAM J. Math. Anal. \textbf{44} (2012), 1894--1924.

\bibitem{NguyenVogelius2}
\bysame, \emph{{Full Range Scattering Estimates and their Application to
  Cloaking}}, Arch. Rational Mech. Anal. \textbf{203} (2012), 769--807.

\bibitem{PayneWeinberge}
L.~E. Payne and H.~F. Weinberger, \emph{{New bounds for solutions of second
  order elliptic partial differential equations}}, Pacific J. Math. \textbf{8}
  (1958), 551--573.

\bibitem{PendrySchurigSmith}
J.~B. Pendry, D.~Schurig, and D.~R. Smith, \emph{{Controlling electromagnetic
  fields}}, Science \textbf{312} (2006), 1780--1782.

\bibitem{Rellich}
F.~Rellich, \emph{{Darstellung der Eigenwerte von $\Delta u+ \lambda u=0$ durch
  ein Randintegral}}, Math. Z. \textbf{46} (1940), 635--636.

\bibitem{RuanYanNeffQiu}
Z.~Ruan, M.~Yan, C.~M. Neff, and M.~Qiu, \emph{{Ideal cylindrical cloak:
  Perfect but sensitive to tiny perturbations}}, Phys. Rev. Lett. \textbf{99}
  (2007), 113903.

\bibitem{Schurig06}
D.~Schurig, J.~J. Mock, J.~Justice, S.~A. Cummer, J.~B. Pendry, A.~F. Starr,
  and D.~R. Smith, \emph{{Metamaterial electromagnetic cloak at microwave
  frequencies}}, Science \textbf{314} (2006), 1133628.

\bibitem{Tricarico-Alu}
S.~Tricarico, F.~Bilotti, A.~Alu, and L.~Vegni, \emph{{Plasmonic cloaking for
  irregular objects with anisotropic scattering properties}}, Phys. Rev. E
  \textbf{81} (2010), 026602.

\bibitem{VasquezMiltonOnofrei}
F.~G. Vasquez, G.~W. Milton, and D.~Onofrei, \emph{{Active Exterior Cloaking
  for the 2D Laplace and Helmholtz Equations}}, Phys. Rev. Lett. \textbf{103}
  (2009), 073901.

\bibitem{VasquezMiltonOnofrei2}
\bysame, \emph{{Broadband exterior cloaking}}, Optics Express \textbf{17}
  (2009), 14800.

\bibitem{WederRigorous}
R.~Weder, \emph{{A rigorous analysis of high-order electromagnetic invisibility
  cloaks}}, J. Phys. A: Math. Theor. \textbf{41} (2008), 065207.

\bibitem{WederRigousTimeDomain}
\bysame, \emph{{The boundary conditions for point transformed electromagnetic
  invisibility cloaks}}, J. Phys. A: Math. Theor. \textbf{41} (2008), 415401.

\bibitem{YanRuanQiu07}
M.~Yan, Z.~Ruan, and M.~Qiu, \emph{{Cylindrical invisibility cloak with
  simplified material parameters is inherently visible}}, Phys. Rev. Lett.
  \textbf{99} (2007), 233901.

\end{thebibliography}

\end{document}